\documentclass{amsart}
\usepackage[foot]{amsaddr}

\title{The Power of an Adversary in Glauber Dynamics}

\author{Byron Chin, Ankur Moitra, Elchanan Mossel, Colin Sandon}
\address{Department of Mathematics, Massachusetts Institute of Technology}
\email{\{byronc, moitra, elmos, csandon\}@mit.edu}

\usepackage{arxiv}

\begin{document}

\maketitle

\begin{abstract}%
Glauber dynamics are a natural model of dynamics of dependent systems. While originally introduced in statistical physics, they have found important applications in the study of social networks, computer vision and other domains. 
In this work, we introduce a model of corrupted Glauber dynamics whereby instead of updating according to the prescribed conditional probabilities, some of the vertices and their updates are controlled by an adversary. We study the effect of such corruptions on global features of the system. Among the questions we study are: How many nodes need to be controlled in order to change the average statistics of the system in polynomial time? And how many nodes are needed to obstruct approximate convergence of the dynamics?

Our results can be viewed as studying the robustness of classical sampling methods and are thus related to robust inference. The proofs connect to classical theory of Glauber dynamics from statistical physics. 
\end{abstract}


\section{Introduction}
Many models of great interest in mathematics have been motivated by the goal of explaining complex systems and their behavior. The phase transitions of magnets, the evolution of opinions on a social network, the spread of infections among a population, and the inter-dependencies between nodes of a biological network are important examples of systems of interest. 

Statistical physics studies the macroscopic behaviors of large systems that emerge. Often times these systems are built from simple building blocks. Major efforts have been invested to describe the properties of such systems in equilibrium, and how they change under perturbations. For example, what is the effect of imposing boundary conditions? And how do small differences in the initialization of the system play out? The seminal works of Martinelli and Toninelli~\cite{martinelli1999lectures, martinelli2010mixing} study the relaxation to equilibrium with various boundary conditions imposed on the Ising model. The study of differences in initialization is often called ``damage spreading" as coined by Kauffman~\cite{kauffman1969metabolic}. Many works have performed physical experiments and simulations to estimate the impact of small changes in initialization on macroscopic properties in equilibrium, see e.g.~\cite{stanley1987dynamics}.  See~\cite{puzzo2008damage} for a survey of the work and results in this area. 

But what can we say about the effect of perturbations over short periods of time as the system evolves toward equilibrium? A natural question to ask in this framework is: what macroscopic behavior can an adversary \textit{efficiently} induce in a physical system with control over some subset of its microscopic constituents? By the term efficiently we mean that we are interested in the effects after at most a polynomial number of time steps of running the natural dynamics. We emphasize that our notion of an adversary is new and flexible, as it allows the nodes to act, collude, and change over time arbitrarily. 

In this work, we focus on systems represented by an Ising model. The Ising model has wide-ranging applications, including as a model for opinion dynamics in a social network~\cite{liu2010influence}, as a model for computer networks~\cite{antonio2010ising}, and as a model for a biological system~\cite{Majewski2001TheIM}. See~\cite{lipowski2022ising} for an extensive collection of works studying the Ising Model from various perspectives. Analogous questions are similarly interesting in other settings such as the Potts model, other Markov random fields and spin glasses. In such settings the macroscopic behavior can be much more complex, and the effectiveness of the adversary becomes even more apparent.

The Glauber dynamics introduced by Glauber~\cite{glauber1963} is the natural formalism by which to add a time component to the Ising model to describe its evolution and progress towards equilibrium. Indeed, prior work by Montanari and Saberi~\cite{montanari2010spread} and Lara et al.~\cite{lara2019analogy} used it to study the spread of innovation and disease across a network. More recently Baldassarri et al.~\cite{baldassarri2022ising} studied the evolution of opinions under Glauber-type dynamics. Another interesting application is \cite{jiang2012fast}, where the Glauber dynamics are used to produce a scheduling algorithm on a wireless network. Each of these examples can be subject to adversarial behavior of collections of nodes or agents.

One setting to highlight is that of a social network. Each vertex of a graph represents a person, with the corresponding spin being their opinion on a specified topic. The edges between vertices correspond to personal connections that may influence the opinions of both people. Encoding opinions in this way allows for the analysis of both the distribution of opinions over the entire population, as well as how it changes over time. The adversary in this case can be interpreted as an advertisement campaign or an external party seeking to influence the global population as efficiently as possible.

\subsection{Our Contributions}
We introduce a general model of adversarially perturbed dynamics called \textit{corrupted Glauber dynamics}. And we establish the first baselines, in this natural model, of the power of adversarial corruptions. 

For context, in the high temperature regime, it is known that the presence of boundary conditions does not alter the equilibrium measure by much. For example, the notion of strong spatial mixing introduced by Martinelli and Olivieri~\cite{martinelli1994approach} stipulates that in equilibrium, the correlation between spins decreases exponentially in the graph distance between the corresponding vertex sets. While in physical systems it is natural to assume that the boundary conditions are ``far away", in our setting, the set of controlled vertices can be distributed in such a way that every node is close to one or more of the corrupted vertices. 
Our goal is to study the effect of controlling general subsets and we show that controlling a small set of corrupted vertices has a small effect on the system as a whole. 

Second, in the low temperature regime, it was proven in many settings, starting with Peierls~\cite{peierls1936ising} that the equilibrium measure consists of two distinct phases, a positive phase and a negative phase. Moreover, we know that the intermediate configurations between these phases are super-polynomially unlikely under the equilibrium measure for many families of graphs, see for example \cite{levin2017markov}. In this strong correlation regime, fixing even a single vertex to $+1$ can heavily bias the configuration toward the positive phase at equilibrium. Surprisingly, this property does not always carryover when we restrict to a polynomial number of steps of the dynamics. In particular, we establish that control of a small constant fraction of the vertices is insufficient to flip the phase in expander graphs.  On graphs that are less expanding such as the $n \times n$ grid, we establish results showing to change the phase in polynomial time requires $\Theta(n)$ vertices. 

Finally, we move away from the setting of fixed boundary conditions 
and take advantage of the flexibility of our definition to accommodate more complex goal. In particular, suppose the corrupted set of vertices behaves adversarially in more sophisticated and coordinated ways. In the low temperature regime, it is known that configurations away from one of the two equilibrium phases are super-polynomially unlikely.  We consider how much control the adversary requires to force the dynamic configuration to remain in these unlikely states, which are most polarized, for super-polynomial time. We prove that in the mean field setting with $n$ vertices, it suffices to control slightly more than $\sqrt{n}$ vertices to induce polarization in the spin configuration for a super-polynomial amount of time. Thus a coordinated adversary can force a system to remain in a fragile state for an extended period of time. 

\subsection{Main Results}
The main purpose of our paper will be to analyze a new process, which we call the corrupted Glauber dynamics, in which some vertices update adversarially and in a coordinated fashion.  

\begin{definition}[Corrupted Glauber Dynamics]\label{def: Corrupted Glauber}
Consider a graph $G$ and a subset of vertices $A \subset V(G)$. The $A$-corrupted Glauber dynamics is a Markov Chain that follows the transitions of the standard Glauber dynamics, except when a vertex $v \in A$ is chosen to be updated, the update is over-ridden by a rule specified by an adversary. The rule can be any function of the current and past states, and has no requirement of computational efficiency. The adversary can update any of the corrupted vertices at any time. We use $m$-corrupted Glauber dynamics to refer to $A$-corrupted Glauber dynamics where $A$ is unspecified other than $|A| = m$. 
\end{definition}

\begin{remark}
In the case of a particular set of vertices being pinned to a fixed state by the adversary, the resulting measure can be interpreted as the Ising model with boundary conditions imposed. See~\cite{martinelli2010mixing,DSVW:04,Weitz:05} for a definition and discussion of boundary conditions.
\end{remark}

The magnetization is the canonical measure of the global state of an Ising model as it measures the average behavior of the system. In much of the discussion below we will assume that the adversary wishes to shift the magnetization in a certain direction. 

\begin{definition}[Magnetization]\label{def: magnetization}
    Given a configuration $\sigma \in \Omega$ the magnetization $m(\sigma) := \sum_{v \in V(G)} \sigma(v)$ is the sum of all the spins over all the vertices. 
\end{definition}

Our main results can be summarized in the following theorems. 
In the high temperature regime, it is known that boundary conditions, i.e. a small set of vertices that are far away from most vertices, induces a negligible change in the resulting equilibrium measure as measured by magnetization. One sense in which this is captured is the property of strong spatial mixing and correlation decay.
In our models, it is not natural to consider only far away vertices, so we consider the more general question of the effect of a small vertex set. The following theorem shows that the two configurations can be coupled such that at any time step, the number of vertices that differ is very likely to be small. In particular, controlling a small number of vertices cannot change the overall statistics by my much in the high temperature setting. 

\begin{theorem}
    Let $G$ be an $n$-vertex, $d$-regular graph $G$. For $\beta(d)$ sufficiently small, when controlling $\epsilon n$ of the vertices, at any fixed time the probability that the magnetization changes by $2\epsilon$ is exponentially small in $n$. In particular the trajectory can be coupled with a standard dynamics such that difference in magnetization is at most $2\epsilon$ at a given time with high probability.
\end{theorem}

Note that the theorem above implies both that in polynomial time it is impossible to deviate by more than $2 \epsilon$ and also that the same is true for the stationary measure. 
See Theorem \ref{thm: high temp} for a formal statement. In contrast, in the regime of low temperature, it is known that the equilibrium exhibits a non-zero magnetization \cite{peierls1936ising}. As a result imposing homogeneous boundary conditions (as in \cite{martinelli2010mixing}), or even fixing a single vertex, will select a phase of the equilibrium measure. 
Thus, it is natural to wonder if a small corrupted set can change the phase of the configuration efficiently. On the contrary, we prove that on an expander graph, a small fraction of the vertices does not suffice to dynamically flip the magnetization in a polynomial number of time steps. Thus there are important differences in the efficacy of adversarial corruptions that vary greatly depending on the inherent robustness and error correction properties of the underlying network. 

\begin{theorem}
    Let $G$ be an $n$-vertex $d$-regular $\alpha$-edge expander graph. For $\beta(d, \alpha)$  sufficiently large and $\epsilon(d, \alpha)$ sufficiently small, when controlling $\epsilon n$ vertices the stationary measure still consists of two phases -- one with positive magnetization and one with negative magnetization. When starting from an all $+1$ configuration it takes exponentially many time steps to move to the negative phase.
\end{theorem}

See Theorem \ref{thm: low temp} for a formal statement. For graphs which are not such good expanders, the behavior is likely more subtle. Controlling one vertex still is not sufficient to change the phase in polynomial time.
However, it suffices to control $o(|V|)$ vertices to change the phase in polynomial time.
For example for the $n \times n$ grid we prove the following: 

\begin{theorem}
    Consider the all $+1$ configuration on the $n \times n$ grid. For $\beta = \Omega(\log n)$ one needs to control linear in $n$ vertices to reach the $-1$ configuration in a polynomial number of steps. On the other hand, for $\beta$ a sufficiently large constant, controlling $o(n)$ vertices yields two phases of the stationary measure, and an exponential number of steps is needed to reach the negative phase.
\end{theorem}

See Theorems \ref{thm: grid infinite}, \ref{thm: grid growing}, and \ref{thm: grid large} for the precise results. It is an interesting question to develop a more refined understanding of the interaction between adversarial corruptions and network topology. But our results already set the stage by establishing some important phenomena. 

Next we ask the question: Are there macroscopic changes that a dynamic adversary, who does not merely pin all the corrupted vertices, can induce? We study the power of an adversary to prolong polarization. We use the term polarization to refer to the state of having an approximately equal number of $+1$ and $-1$ spins, in other words a magnetization close to zero, which is typically an exponentially unlikely configuration. The following theorem gives a tight bound on the size of the corrupted set necessary for the adversary to maintain polarization for longer than polynomial time.

\begin{theorem}\label{thm: convergence informal}
    For $\beta$ a sufficiently large constant in the mean field setting, starting from an i.i.d. random configuration, controlling $n^{0.5+\epsilon}$ vertices maintains polarization for super-polynomial (in other words, it prevents convergence to one of the phases in polynomial time). 
\end{theorem}

To prove this theorem, a more sophisticated and coordinated attack is needed. Our attack is based on assigning the corrupted vertices to the minority spin, until that spin nearly becomes an overwhelming majority. This is formalized in Theorem \ref{thm: convergence}. Note that the random initialization in this result is essential, as corrupting $o(n)$ vertices starting from a global $+1$ configuration has no desired effect. Understanding the scope and extent of an adversary's power to shape opinion dynamics is obviously a challenging problem. Our results show that, depending on the goal, an adversary's strategy can be subtle and hard to detect. 

Another interesting question is to study our model from the perspective of the adversary, namely the optimization problem of maximizing its effect. For example, for ferromagnetic Ising models, given a budget of $\epsilon n$ nodes, which ones should the adversary choose to minimize the magnetization when measured at a specific time in the future? This question is closely related to the question of influence 
maximization~\cite{KeKlTa:03,MosselRoch:10}. We note that for the ferromagnetic Ising model if we are interested in the long term-behavior (stationary) behavior of the system then it is known that the influence maximization problem is sub-modular, see~\cite{BrKoMo:19}. Our main result in this direction is that the sub-modularity extends to the dynamic setting under a sufficiently strong external field. 

\begin{theorem}\label{thm: dynamic submodular informal}
    For the corrupted Glauber dynamics on a graph $G$ with external field at least the maximum degree, the influence maximization problem at any time step $t$ is submodular. 
\end{theorem}

This is formalized in Theorem \ref{thm: dynamic submodular}. The proof captures a general phenomenon that submodularity of the local update rule implies global submodularity of the maximization problem. 

\subsection{Future Questions}
Our framework leaves many interesting questions and research directions.

Is it possible to generalize our result about preventing convergence outside of the mean field setting? It is natural to believe that the conclusion should hold true as long as the graph is sufficiently well-connected, which leads to the following conjecture. 

\begin{conjecture}
Let $G$ be an expander graph. Then there is a choice of corrupted set such that the conclusion to Theorem \ref{thm: convergence informal} holds. 
\end{conjecture}

We expect many interesting phenomena to occur for models that allow anti-ferromagnetic interactions, including spin glasses, where we believe that finding the optimal set (even approximately) is computationally hard. Indeed given that sampling from such models is computationally hard~\cite{Sly:10,SlySun:12}, it is natural to expect that influence maximization is also hard. 

Another dimension of complexity may be introduced by considering game theoretic models where a number of different players with different objectives may choose different sets. For example, what are the optimal strategies for $3$ players in an anti-ferromagnetic Potts models if one wants to maximize the number of green nodes, another the number of blue nodes and the third the number of red nodes at a specific time? We can assume for example that players each have a budget of $\epsilon n$ nodes and a node is allocated uniformly at random to one of the players that bid for it. 

It is interesting to answer analogous questions here, where the study of censored rather than corrupted nodes has also been of recent interest, see e.g. \cite{moitra2021learning} and references within.

\subsection*{Acknowledgements}
BC was supported by the NSF Graduate Research Fellowship Program. 
EM and CS were partially supported by Vannevar Bush Faculty Fellowship ONR-N00014-20-1-2826 and by 
Simons Investigator award (622132)
EM was also partially supported by NSF award CCF 1918421. AM and CS were also partially supported by AM's ONR YIP grant.  AM was also partially supported by a David and Lucile Packard Fellowship.

\section{Preliminaries}
Here we introduce the main definitions and notation that we will use for the rest of the paper. 
\begin{definition}[Ising Model]\label{def: Ising model}
    Let $G = (V(G), E(G))$ be a graph and $\Omega = \{+1, -1\}^{V(G)}$. The Ising Model with parameter $\beta$ is defined by the probability measure $\mu$ on $\sigma = (\sigma(v))_{v \in V(G)} \in \Omega$ where
    \[ \mu(\sigma) = \frac{e^{\beta H(\sigma)}}{Z(\beta)} \]
    where $H(\sigma) = \sum_{uv \in E(G)}\sigma(u)\sigma(v)$ and $Z(\beta) = \sum_{\sigma \in \Omega} e^{\beta H(\sigma)}$ is the normalization constant.
\end{definition}

The parameter $\beta$ is known as the inverse temperature, which governs the strength of the dependencies between spins on adjacent vertices in $G$. This distribution will be a vehicle to study and quantify macroscopic changes that can be induced when certain sets of vertices are controlled by an adversary. 

In order to study the process of how the effect of a corrupted set of vertices permeates through the configuration, we need a process for re-sampling the spins of vertices that agrees with the Ising model distribution. This role will be played by the well-known Glauber dynamics.

\begin{definition}[Glauber Dynamics]\label{def: Glauber}
    The Glauber dynamics is a Markov chain $X_t$ on $\Omega$ with the following transition rule. Choose a vertex $v \in V(G)$ uniformly at random. Then $X_{t+1}(u) = X_t(u)$ for every $u \neq v$ and $X_{t+1}(v) = +1$ with probability $\frac{\exp(\beta\sum_{uv \in E(G)} X_t(u))}{\exp(\beta\sum_{uv \in E(G)} X_t(u)) + \exp(-\beta\sum_{uv \in E(G)} X_t(u))}$.
\end{definition}

The transition probability is the conditional probability of the spin at $v$ being $+1$ under the measure $\mu$ conditioned on the remaining spins in the graph. It is easy to check that the Glauber dynamics is reversible with respect to the measure $\mu$, and thus has our desired stationary distribution. 

The Glauber dynamics is one of the most extensively studied ways to sample from the Ising model distribution. Given the computational hardness of exactly sampling from the measure (see e.g. \cite{istrail2000statistical}), such a Markov chain is the most common way to sample approximately. There is a rich line of research aimed at analyzing the Glauber dynamics. Through many innovations and breakthroughs, it is known that in many cases, the chain exhibits a phase transition in the time it takes to reach equilibrium. For a comprehensive overview of techniques see \cite{martinelli1999lectures} and \cite{levin2017markov}. 

At a technical level, some of our results will be established using bottlenecks for the magnetization. 
More generally, a bottleneck is a constraining region of the configuration space $\Omega$ that limits the transitions of the Glauber dynamics. The strength of a bottleneck is quantified by the bottleneck ratio. 

\begin{definition}[Bottleneck Ratio]\label{def: bottleneck}
    Let $P$ be the transition matrix of the Glauber dynamics on $\Omega$, which has stationary distribution $\mu$. Define 
    \[ Q(\sigma_1,\sigma_2) := \mu(\sigma_1)P(\sigma_1,\sigma_2) \text{ and } Q(A, B) = \sum_{\sigma_1\in A, \sigma_2\in B} Q(\sigma_1,\sigma_2). \]
    The bottleneck ratio of a subset $S \subset \Omega$ is defined as 
    \[ \Phi(S) := \frac{Q(S, S^c)}{\mu(S)} \]
    and the bottleneck ratio of the dynamics is 
    \[ \Phi^* := \min_{\mu(S) \leq \frac{1}{2}} \Phi(S). \]
\end{definition} 

As alluded above, the bottleneck ratio controls the time it takes for the Glauber dynamics to move from one region of $\Omega$ to another, and will be used to bound the number of steps needed for the impact of the corrupted set $A$ to take effect. See \cite{levin2017markov} Theorem 7.4 for a useful result in this spirit. 

\subsection*{Notation} 
For the rest of the paper, we will use the convention that variables with a tilde will refer to the corrupted dynamics (e.g. $\tilde{X_t}, \tilde\mu, \tilde\Phi$, etc) and variables with no tilde will refer to the corresponding values in the standard dynamics. 

\section{No Corruption Effect}
We first prove ``negative" results (in the sense that the corrupted vertices do not change much) in both the high and low temperature regimes. These results will all take place on regular graphs. The first result is in the case of a high temperature (small $\beta$) system, where interactions between vertices are close to independent. In this case it is natural to believe that a group of few vertices cannot strongly influence the remaining vertices. This is captured in a quantitative sense by comparing the trajectories between $X_t$ and $\tilde X_t$.
\begin{theorem}\label{thm: high temp}
    Fix $\epsilon > 0$. Suppose $d\tanh\beta \ll 1$, and let $G$ be a $d$-regular graph. Let $X_t$ be the standard Glauber dynamics, and $\tilde X_t$ be $\epsilon n$-corrupted Glauber dynamics coupled in the natural way. Let $D_t = \{v \in V(G) \setminus A: X_t(v) \neq \tilde X_t(v)\}$. Then for any time $t \geq 0$, 
    \[ \Prob{|D_t| > 2 \epsilon n} \leq 8 \beta\epsilon e^{-2\epsilon^2\beta n}.\]
\end{theorem}

The proof proceeds by comparing the trajectories of the standard dynamics with the corrupted one. We show that when the number of differences is large, there is a bias toward returning to 0. Thus, at any fixed time, it will be unlikely for there to be a large number of differences between the two processes. The first claim is captured by the following lemma.

\begin{lemma}\label{lemma: drift}
For any $\delta > 0$ let $L_1 = n\left( \delta + 
\frac{(\epsilon + \delta)d\beta}{1-d\beta} \right)$. Then for all $L \ge L_1$, 
\[ \E{|D_{t+1}| \middle\vert |D_t| = L} \leq L - \delta.\] 
\end{lemma}

\begin{proof}
Suppose that $|D_t| = L$. In order for $D_t$ to decrease, we need the Glauber dynamics to select a vertex in the difference set. This happens with probability $\frac{L}{n}$. Then to decrease, the updates must match, which occurs with probability $1 - \abs{\mu_{X_t}(\sigma_v=+) - \mu_{\tilde X_t}(\sigma_v=+)}$. 
In order for $D_t$ to increase, it needs to select a vertex outside the difference set with probability $1-\epsilon - \frac{L}{n}$, and produce updates that do not match with probability $\abs{\mu_{X_t}(\sigma_v=+) - \mu_{\tilde X_t}(\sigma_v=+)}$. 
The remaining $\epsilon$ probability is when Glauber selects a vertex in the corrupted set, at which $D_t$ does not change. 

To improve the dependence on $d$ and $\beta$ and obtain the desired $L_1$, we further distinguish the cases by the number of differences in the neighborhood of the selected vertex $v$ to compute the above probabilities. Let $$V_t(k) = \{v: N_{X_t}(v) \text{ and } N_{\tilde X_t}(v) \text{ contain $k$ differences}\}.$$ By double counting, \[ \sum_{k \geq 0} k \cdot |V_t(k)| \leq (\epsilon n + L)d. \]
Now, if the neighborhood of $v$ has $k$ differences, the difference in probabilities is $$\frac{1}{2}\tanh(\beta(s+k)) - \frac{1}{2}\tanh(\beta(s-k)) \leq \tanh(k\beta).$$ Thus we want 
\begin{align*}
    &\frac{1}{n} \sum_{k \geq 0} \abs{D_t \cap V_t(k)} \cdot (1 - \tanh(k\beta)) - \abs{D_t^c \cap V_t(k)} \cdot \tanh(k\beta) \\
    = \quad & \frac{1}{n} \sum_{k \geq 0} \abs{D_t \cap V_t(k)} - \abs{V_t(k)} \tanh(k\beta) > \delta.
\end{align*}
Since $\tanh(x) < x$ for $x > 0$, it suffices for
\[ \frac{L}{n} - \frac{1}{n}\sum_{k \geq 0} \abs{V_t(k)} \cdot k\beta 
    \geq \frac{L}{n} - \frac{\beta(\epsilon n + L)d}{n} 
    = \frac{L(1-d\beta)}{n} - \epsilon d\beta > \delta. \]
This holds when
\[ L \geq \frac{n(\delta + \epsilon d\beta)}{1 - d\beta} = n\left( \delta + 
\frac{(\epsilon + \delta)d\beta}{1-d\beta} \right). \]
\end{proof}

The following lemma is an estimate on the hitting probabilities of a biased random walk.

\begin{lemma}\label{lemma: biased walk}
    Fix $0 < \delta < \frac{1}{4}$. Let $W_t$ be a random walk such that $W_0 = 1$, $\Prob{W_{t+1} = W_t + 1} = \frac{1}{2}-\delta$. Let $\tau_a = \inf\{t \geq 0: W_t = a\}$. Then 
    \[ \Prob{\tau_a < \tau_0} < 8\delta e^{-2\delta a}. \]
\end{lemma}

\begin{proof}
By considering the martingale $M_t = \left( \frac{\frac{1}{2} + \delta}{\frac{1}{2} - \delta} \right)^{W_t}$ and the stopping time $\tau = \min\{\tau_a, \tau_0\}$ we can compute
\begin{align*}
    \frac{\frac{1}{2} + \delta}{\frac{1}{2} - \delta} = \E{M_0} = \E{M_\tau} = \Prob{\tau_0 < \tau_a} + \Prob{\tau_a < \tau_0}\cdot \left( \frac{\frac{1}{2} + \delta}{\frac{1}{2} - \delta} \right)^a
\end{align*}
Rearranging gives 
\[ \Prob{\tau_a < \tau_0} = \frac{\frac{\frac{1}{2} + \delta}{\frac{1}{2} - \delta} - 1}{\left( \frac{\frac{1}{2} + \delta}{\frac{1}{2} - \delta} \right)^a - 1} < \frac{8\delta}{(1 + 4\delta)^a -1} < 8\delta e^{-2\delta a}. \]
\end{proof}

\begin{proof}[Proof of \cref{thm: high temp}] 
Set $\delta = \beta\epsilon$ and let $L_1$ be as in Lemma \ref{lemma: drift}. Define $T_0 = \inf\{t: |D_t| > L_1\}$ and $S_0 = \inf\{t > T_0: |D_t| = L_1\}$. Inductively define $T_k = \inf\{t > S_{k-1}: |D_t| > L_1\}$ and $S_k = \inf\{t > T_k: |D_t| = L_1\}$, so that the time intervals $[T_k, S_k)$ capture the excursions of $|D_t|$ above $L_1$. For each excursion, associate a new process $(W_t^{(k)})_{t=0}^{S_k-T_k}$ started from 0 that increases by 1 with probability $\frac{1}{2} - \delta$ and decreases by 1 with probability $\frac{1}{2} + \delta$ whenever $|D_t|$ moves. By Lemma \ref{lemma: drift}, we can couple the excursion $(|D_t|)_{t=T_k}^{S_k}$ with $W_t^{(k)}$ such that $|D_t| \leq W_{t-T_k} + L_1$. Since each $W_t^{(k)}$ is an identically distributed biased random walk we know that for each $k$, by Lemma \ref{lemma: biased walk},
\[ \Prob{\exists t\leq S_k-T_k: W_t^{(k)} \geq \epsilon n} < 8\delta e^{-2\delta \epsilon n}. \]
Now, for any fixed time $t$, if $|D_t| \geq L_1 + \epsilon n$, then by the invariant above it must be coupled to a biased walk satisfying $W_{t'}^{(k)} \geq \epsilon n$ for the corresponding $k$ and $t'$. This implies that 
\[ \Prob{|D_t| \geq 2\epsilon n} < \Prob{|D_t| \geq L_1 + \epsilon n} \leq \Prob{W_{t'}^{(k)} \geq \epsilon n} < 8\epsilon\beta e^{-2\epsilon^2\beta n} \]
which is the desired inequality. 
\end{proof}

We now consider the opposite case of a low temperature (large $\beta$) system. In this case, there are strong interactions between adjacent vertices, so it is possible to hope for a large impact on the dynamics while only controlling a small number of vertices. However, as long as the graph is sufficiently well-connected (captured by the notion of edge expansion), then a small constant fraction of the vertices is still insufficient to flip the phase of the configuration efficiently. This is captured by the following theorem. 

\begin{theorem}\label{thm: low temp} 
    Let $G$ be a $d$-regular, $\alpha$-edge-expander graph. There exists a $C(\frac{\alpha}{d})$ such that for $\beta > \frac{1}{2d} C(\frac{\alpha}{d})$ and $\epsilon < \frac{\alpha}{4d}$, the mixing time of the $\epsilon n$-corrupted Glauber dynamics is at least $C\exp(\frac{1}{4}\alpha\beta n)$. In particular, if $\alpha = \Theta(d)$ then this holds for $\beta = \Omega(\frac{1}{d})$. 
\end{theorem}

\begin{proof}
The proof proceeds by exhibiting a bottleneck in the dynamics of the process. Suppose $G$ a $d$-regular expander with $d \geq \min_{0 \leq |S| \leq \frac{n}{2}} \frac{|\partial S|}{|S|} \geq \alpha$. Define $A_k = \{ \sigma: \abs{\{v: \sigma_v = +\}} = k\}$. For $\sigma \in A_k$, we know $dk \geq \abs{\partial S} \geq \alpha k$ where $S = \{v: \sigma_v = +\}$. Note that 
\[ \sum_{u \sim v} \sigma_u\sigma_v = (|E(G)| - |\partial S|) - |\partial S| \in \left[ \frac{nd}{2} - 2dk, \frac{nd}{2} - 2\alpha k \right]. \]

Now, consider a dynamics that is corrupted with $\epsilon n$ vertices fixed to $+$. Here, with only $n(1-\epsilon)$ free vertices, we adjust the bottleneck to $k^* = \frac{n(1-2\epsilon)}{2}$. Note that the number of configurations with $k < k^*$ is less than that of $k > k^*$, and also the configurations with $k < k^*$ have lower probabilities due to the fixed set of $+$ vertices. Thus, $\tilde\mu\left(\bigcup_{k < k^*} A_k\right) < \frac{1}{2}$. With $k+\epsilon n$ vertices labeled $+'s$, 
\begin{align*}
    \tilde\mu(A_{k^*}) &\leq \frac{1}{Z_c} \cdot \binom{n(1-\epsilon)}{\frac{n(1-2\epsilon)}{2}}\exp(\beta n\left(\frac{d}{2} - \alpha \right)) =: \frac{1}{Z_c} \cdot \tilde{a}_{k^*} \\
    \tilde\mu\left(\bigcup_{k < k^*} A_k\right) &\geq \max_{k > k^*} \frac{1}{Z_c} \cdot \binom{n(1-\epsilon)}{k}\exp(\beta\left(\frac{nd}{2} - 2d(k + \epsilon n)\right)) =: \max_{k < k^*} \frac{1}{Z_c} \cdot \tilde{a}_k
\end{align*}
By the same asymptotic approximations, and writing $k = cn(1-\epsilon)$,
\begin{align*}
    \log \tilde{a}_{cn(1-\epsilon)} &= \beta nd\left(\frac{1}{2} - 2(c(1-\epsilon) + \epsilon)\right) - n(1-\epsilon)[c\log c + (1-c)\log (1-c) + o(1)] \\
    \log \tilde{a}_{k^*} &\leq n(1-\epsilon) \log 2 + \beta n(\frac{d}{2} - \alpha)
\end{align*}
The first expression has the same critical point at $c^* = \frac{1}{1+e^{2\beta d}}$, so we can bound the log of the bottleneck ratio as 
\begin{align*}
    \log \tilde{\Phi}^* &\leq n\left[ \beta(\frac{d}{2} - \alpha) + (1-\epsilon)\log 2\right] \\
    &\qquad - n\left[ \beta d(\frac{1}{2} - 2(c^*(1-\epsilon) + \epsilon)) - (1-\epsilon) [c^*\log c^* + (1-c^*)\log (1-c^*)] \right] \\
    &\leq n[(1-\epsilon)(\log 2 + c^*\log c^* + (1-c^*)\log (1-c^*)) + \beta(-\alpha + 2d(c^*(1-\epsilon) + \epsilon))].
\end{align*}
The rightmost term is at most $-\frac{1}{2} n\alpha\beta$ as long as $c^* < \frac{\alpha}{4d(1-\epsilon)}$ and $\epsilon < \frac{\alpha}{4d}$. The first condition translates to $\beta > \frac{1}{2d} \log(\frac{8d}{\alpha} - 1)$.

To handle the remaining term, we set $C := 2d\beta$. We want to ensure that 
\[ \log 2 + c^*\log c^* + (1-c^*)\log (1-c^*) \leq \frac{\alpha \beta}{4} = \frac{C\alpha}{8d}. \]
Taking the derivative of the left hand side with respect to $C$, we get the expression $\frac{Ce^C}{(1+e^C)^2}$. Notice that as $C \to \infty$, this derivative tends to zero, while $\frac{d}{dC} \frac{C\alpha}{8d} = \frac{\alpha}{8d}$. Thus, there exists some $C(\frac{\alpha}{d})$ such that for all $C > C(\frac{\alpha}{d})$, the desired inequality holds. In particular, for $\beta \geq \frac{C(\alpha/d)}{2d}$ we have $\log \tilde \Phi^* \leq -\frac{\alpha\beta}{4}$. 

This gives the desired bottleneck in the corrupted case. The result for mixing time follows from a standard inequality (i.e. \cite{levin2017markov} Theorem 7.4).
\end{proof}

\section{Forcing on the Grid}
In this section we give a series of ``positive" results in which the corruption has a macroscopic effect the dynamics in polynomial time. These results will take place on the $n \times n$ square lattice. As opposed to the prior results which are for general graphs satisfying certain nice properties, the analysis of Glauber dynamics on the grid rely on the specific geometric structure of the graph. The common question we study is: starting from the $+1$ configuration, how many nodes need to be corrupted in order to force the dynamics to the $-1$ configuration in a polynomial (in $n$) number of steps? 

\begin{remark}
By the monotonicity of the measure $\mu$, we only need to consider the case in which the corrupted vertices are fixed to the $-1$ spin. 
\end{remark}

To start, we consider the case of zero temperature, which is closely tied to the study of bootstrap percolation (see \cite{balogh1998random} for further background on this related problem). 

\begin{theorem}\label{thm: grid infinite}
Let $\beta = \infty$ and $G$ be the square lattice with vertex set $[n]^2$ . Suppose the initialize with the identically $+1$ configuration. Then:
\begin{enumerate}
    \item The $m$-corrupted Glauber dynamics can never reach the $-1$ configuration for $m < n$.
    \item Let $A = \{(i,i): i \in [n]\}$. Then the $A$-corrupted Glauber dynamics will reach the $-1$ configuration in expected time $O(n^5)$. 
\end{enumerate}
\end{theorem}

\begin{proof}
Note that at zero temperature, the vertices that can eventually transition to a $-1$ state is exactly the bootstrap percolation set of the corrupted vertices. Thus, the first claim follows from a perimeter argument (\cite{balogh1998random} Folklore Fact). 

For the second claim, we relate the progress of the configuration toward the $-1$ configuration to a random walk. Let $C_j = \{(i,j): i < j\}$ and $T_j = \bigcup_{k \leq j} C_k$. Notice that once the configuration is identically $-1$ on $T_j$, then this portion of the configuration is frozen in this state. Define $S_j = \inf\{t \geq 0: T_j \text{ is identically $-1$ at time $t$}\}$. We bound $S_{j+1} - S_j$ by comparison to a random walk. Notice that $C_{j+1}$ must flip to the $-$ state in the order $(j, j+1), (j-1, j+1), \ldots, (1, j+1)$. Moreover, if $(j, j+1), \ldots, (k, j+1)$ are all in the $-1$ state, then $(j, j+1), \ldots, (k+1, j+1)$ are all frozen in the $-1$ state. Thus, the quantity $X_t = \inf\{k: (k, j+1) \text{ is $-$ in $\sigma_t$}\}$ is a simple random walk on $[j]$ reflected at $j$ and absorbed at $1$. By the Glauber dynamic update, $\Prob{X_{t+1} = X_t + 1} = \Prob{X_{t+1} = X_t - 1} = \frac{1}{2n^2}$. Starting from $j$ the expected time to absorption of $X_t$ is $j^2$, so $\E{S_{j+1} - S_j} = 2n^2j^2$. Thus, 
\[ \E{S_n} = \sum_{j=1}^{n-1} \E{S_{j+1} - S_j} = \sum_{j=1}^{n-1} 2n^2j^2 = O(n^5). \]
By symmetry, the same argument applies to the lower triangle below the diagonal, and the hitting time of the $-1$ configuration is the maximum of the upper triangle and lower triangle, yielding the desired bound. 
\end{proof}

As a corollary of the above proposition, we can study the case of slowly decreasing temperature. The main observation is that at sufficiently low temperature, the behavior of the chain is unlikely to differ from that of a zero temperature chain during a polynomial number of steps. 

\begin{theorem}\label{thm: grid growing}
Let $\beta \geq 4\log n$ and $G$ be the square lattice with vertex set $[n]^2$ . Suppose the initialize with the identically $+1$ configuration. Then for $A = \{(i,i): i \in [n]\}$, the $A$-corrupted Glauber dynamics will reach the $-1$ configuration in expected time $O(n^5)$. 
\end{theorem}

\begin{proof}
The conditional probability of a vertex whose neighborhood has sum $s > 0$ to be updated to a $-$ state is $\frac{e^{-\beta s}}{e^{\beta s} + e^{-\beta s}} = \frac{n^{-s}}{n^s + n^{-s}} = \frac{1}{n^{8s} + 1} < \frac{1}{n^8}$. By symmetry the conditional probability of a vertex whose neighborhood has sum $s < 0$ to be updated to a $+1$ state is also $< \frac{1}{n^8}$. Denote these events as ``deviations" from the $\beta = \infty$ case studied above in Proposition \ref{thm: grid infinite}. Let $T$ be the first time that the configuration is identically $-1$ in the setting of Proposition \ref{thm: grid infinite}, and let $N$ be the time of the first deviation in the process. Then $N$ dominates $N'$, a $\mathrm{Geom}(n^{-8})$ random variable that is independent of T. We have the estimate 
\[ \Prob{N < n^6} \leq \Prob{N' < n^6} = 1 - \left(1-\frac{1}{n^8}\right)^{n^6} < n^{-2}. \]
We can bound
\begin{align*}
    \Prob{N' < T} &= \sum_{m < n^6} \Prob{m < T}\Prob{N' = m} + \sum_{m \geq n^6} \Prob{m < T} \Prob{N' = m} \\
    &\leq \Prob{N < n^6} + \sum_{m \geq n^6} \Prob{n^6 < T} \Prob{N' = m} \\
    &\leq n^{-2} + \frac{\E{T}}{n^6}\cdot \Prob{N' \geq n^6} \\
    &\leq n^{-2} + Cn^{-1}
\end{align*}
where in the last step we used Markov's Inequality and Proposition \ref{thm: grid infinite}. Thus, with probability $1-o(1)$, we do not deviate before reaching the $-1$ configuration. In the event that a deviation occurs, we ``reset", and by monotonicity the expected time to reach the $-1$ configuration from this new configuration must also be bounded by $O(n^5)$. Conditioning on the number of runs until we reach the $-1$ state without a deviation, we obtain the claimed result.  
\end{proof}
\begin{remark}
    It follows from the proof that once the all $-1$ state is reached, the configuration remains in this state for at least $n^8$ time in expectation. In particular, since the time taken to reach the state is $O(n^5)$ on average, the configuration will remain in the $-1$ state for $1-o(1)$ fraction of time.
\end{remark}
\begin{remark}
    It does not follow from the perimeter argument cited previously that $n-1$ vertices does not suffice for a polynomial time phase flip. Indeed, given $n-1$ vertices as the set $A$, a deviation at the remaining diagonal site will occur in expected polynomial time, after which the global spin configuration will flip in polynomial time once again.
\end{remark}

Finally, to handle the case of low (but constant) temperature, considerably more work is needed. In this case, we do not hope to converge to the $-1$ configuration, but rather to the $-1$ phase, which is the distribution on configurations under $-1$ boundary conditions. The analysis in this case builds on a series of results by Martinelli which analyze the standard Glauber dynamics on the square grid under various boundary conditions. In this setting we show a transition in behavior when the size of the corrupted set moves from sub-linear to linear. The proofs of these results rely on the geometric structure of the lattice. 

\begin{theorem}\label{thm: grid large}
Let $G$ be the $n \times n$ grid, and let $\beta$ be a sufficiently large constant. Starting from the all $+1$ configuration, the $4n$-corrupted Glauber dynamics vertices converges to the $-1$ phase in sub-exponential time. Conversely, the $o(n)$-corrupted Glauber dynamics mixes in exponential time. 
\end{theorem}

\begin{proof}
As shown in \cite{martinelli2010mixing} imposing $-1$ boundary conditions on the $n \times n$ grid yields a mixing time of at most $\exp(cn^\epsilon)$. Thus, the all $+1$ configuration can be forced into the $-1$ phase in sub-exponential time with control of $O(n)$ vertices. 

Conversely, with control of $o(n)$ vertices, the bottleneck that yields exponential time mixing for the square grid remains, so it is not possible to force the $+1$ phase to the $-1$ phase in sub-exponential time. We will make repeated use of the following theorem from \cite{martinelli1999lectures}, Proposition 6.3. 

\begin{theorem*}
For $\beta$ sufficiently large, and $\Delta \in (-m^*(\beta), m^*(\beta))$, 
\[ \lim_{n \to \infty} -\frac{1}{\beta(2n+1)}\log(\mu(\{\sigma: m(\sigma) = \Delta n^2)\}) = c_\Delta > 0 \]
where $m^*(\beta)$ is a constant known as the spontaneous magnetization.  
\end{theorem*}

We use this to construct a bottleneck in the corrupted dynamics. Note that by symmetry, in the standard dynamics with $S = \{\sigma: m(\sigma)>0\}$, we have $\mu(S) = \frac{1 - e^{-c\beta(2n+1)}}{2} < \frac{1}{2}$ and $\Phi(S) < 4e^{-c\beta(2n+1)}$ so the bottleneck ratio is exponentially small in $n$. In the corrupted dynamics, we once again have $\tilde\mu(S) < \frac{1}{2}$.

The idea of the argument is that since the maximum degree of $G$ is 4, the corrupted set of vertices can affect at most $o(n)$ edges, which means the probability of each state $\sigma$ can change by at most a factor of $e^{o(n)}$. Thus, we have that $\tilde\Phi(S) < 4e^{-c\beta(2n+1) + o(n)}$, which also yields an exponential mixing time of at least $e^{c\beta(2n+1) - o(n)}$. This is formalized below to conclude the proof. 

We begin with a lemma that controls the change in the normalization constant of the two distributions. 
\begin{lemma}\label{lemma: Z ratio}
    $Z(\beta) e^{-o(n)} \leq \tilde Z(\beta) \leq Z(\beta).$
\end{lemma}

\begin{proof}
    The upper bound is clear, since the configurations in the corrupted model are a strict subset of the configurations for the standard dynamics, so we prove the lower bound. For $\tau \in \{+1, -1\}^A$, let $\Omega_\tau = \{\sigma: \sigma(A) = \tau\}$. We can re-write
    \begin{align*}
        Z(\beta) &= \sum_\Omega e^{-\beta H(\sigma)} = \sum_\tau \sum_{\Omega_\tau} e^{-\beta H(\sigma)} \leq \sum_\tau e^{4\beta|A|} \tilde Z(\beta) \leq 2^{|A|}e^{4\beta|A|} \tilde Z(\beta) = e^{o(n)} \tilde Z(\beta).
    \end{align*}
    The first inequality follows from the observation that $\tau = -1$ corresponds exactly to $\tilde Z(\beta)$. Since the configurations $\sigma$ and $\sigma_\tau$ obtained by replacing the spins on $A$ with $\tau$ differ in at most $4|A|$ edges, and thus $H(\sigma) - 4|A| \leq H(\sigma_\tau) \leq H(\sigma) + 4|A|$. 
\end{proof}

Let $\tilde m(\sigma)$ be the magnetization on the un-corrupted set of vertices. Define $\tilde S = \{\sigma: \tilde m(\sigma) > 0\}$. Note that this represents half of the possible configurations, and since the corrupted vertices have spins fixed to $-1$, it follows that $\tilde \mu(\tilde S) < \frac{1}{2}$. Note that $\tilde Q(\tilde S, \tilde S^c) \leq \tilde \mu(\{\sigma: \tilde m(\sigma) = 0\})$, so we proceed by bounding the right hand side. Note that $\tilde m(\sigma) = 0$ implies $m(\sigma) = -|A|$ with the assumption of $-1$ spins on the corrupted set. Thus, 
\begin{align*}
\tilde Z(\beta)\tilde \mu(\{\sigma: \tilde m(\sigma) = 0\}) &= \sum_{\tilde m(\sigma) = 0} e^{-\beta H(\sigma)} \leq \sum_{m(\sigma) = -|A|} e^{-\beta H(\sigma)} \\
&= Z(\beta)\mu(\{\sigma: m(\sigma) = -|A|\}) \leq Z(\beta)e^{-c_A\beta(2n+1)}
\end{align*}
for $n$ sufficiently large. It follows from Lemma \ref{lemma: Z ratio} that
\[\tilde \mu(\{\sigma: \tilde m(\sigma) = 0\}) \leq \frac{Z(\beta)}{\tilde Z(\beta)}e^{-c_\Delta\beta(2n+1)} = e^{-c_A\beta(2n+1) + o(n)}. \]
To conclude, it suffices to show that $\tilde \mu(\tilde S)$ is not too small. Note that $m(\sigma) > |A|$ implies $\tilde m(\sigma) > 0$, as the  contribution from the corrupted set to the magnetization is bounded by $-|A|$. Thus,
\begin{align*}
\tilde Z(\beta)\tilde \mu(\tilde S) &= \sum_{\sigma \in \tilde S} e^{-\beta H(\sigma)} \geq \frac{1}{2^{|A|}e^{o(n)}}\sum_{m(\sigma) > |A|} e^{-\beta H(\sigma)} \\
&= \frac{Z(\beta)}{2^{|A|}e^{o(n)}}\mu(\{\sigma: m(\sigma) > |A|\}) \\
&\geq \frac{Z(\beta)}{2^{|A|}e^{o(n)}}\left( \frac{1}{2} - |A|e^{-c_A\beta(2n+1)} \right) \geq Z(\beta) e^{-o(n)}
\end{align*}
Again by Lemma \ref{lemma: Z ratio}, we can deduce that $\tilde\mu(\tilde S) \geq e^{-o(n)}$. Thus, 
\[ \tilde \Phi^* \leq \tilde \Phi(\tilde S) \leq \frac{e^{-c_A\beta(2n+1) + o(n)}}{e^{-o(n)}} = e^{-c_A\beta(2n+1) + o(n)}. \]
\end{proof}

\section{Efficient Polarization}
To conclude we examine a question from a new perspective. As opposed to the above results in which the goal was to influence the overall state of the network as much as possible in a specified direction, we now consider the goal of introducing as much disorder as possible. If the spins of the configuration are initialized randomly, there will be a strong misalignment within the system. We desire the effect of oscillating the configuration repeatedly, perpetuating the polarization of the configuration. In this case, it is also natural to ask how much of the network must be corrupted to induce this behavior. 

\begin{theorem}\label{thm: convergence}
For $G=K_n$, starting from an i.i.d. configuration, the  $n^{0.5 + \epsilon}$-corrupted Glauber dynamics has $m(\tilde X_t)$ fluctuating between $-\sqrt{n}$ and $\sqrt{n}$ up to time $T$ 
with probability at least $(1-o(1))(1 - e^{-\beta n})^{T}$. In particular, for any constant $\beta$ the dynamics will not converge in sub-exponential time.
\end{theorem}

\begin{remark}
    Note that the size of the corrupted set is tight, as a randomly initialized configuration has fluctuations in magnetism of order $\sqrt{n}$, and in order to influence the network the control of the adversary must dominate this fluctuation. 
\end{remark}

\begin{proof}
At initialization, set the spins on $A$ to be identically $+1$. Let $M_t = m(\tilde X_t(V(G) \setminus A))$. When $M_t \geq -\frac{1}{2} n^{0.5+\epsilon}$ we have the following estimates. Suppose that a vertex $v$ is selected by the dynamics and $\sigma(v) = -1$. Then the conditional probability that the spin is updated to $+1$ is (recall Definition \ref{def: Glauber})
\[ \frac{e^{\beta(|A|+M_t+1)}}{e^{\beta(|A|+M_t+1)} + e^{-\beta(|A|+M_t+1)}} = 1 - \frac{1}{e^{2\beta(|A|+M_t+1)} + 1} \geq 1 - e^{-\beta n^{0.5+\epsilon}} =: p. \]
Similarly, if $v$ is selected with $\sigma(v) = +1$, the conditional probability that the spin is updated to $-1$ is 
\[ \frac{e^{-\beta(|A|+M_t-1)}}{e^{\beta(|A|+M_t-1)} + e^{-\beta(|A|+M_t-1)}} = \frac{1}{e^{2\beta(|A|+M_t-1)} + 1} \leq e^{-\beta n^{0.5+\epsilon}} =: q. \]

Thus, in this regime $M_t$ stochastically dominates a random walk with $p$ probability of a $+1$ step and $q$ probability of a $-1$ step. Thus, started from $-n^{0.5+\epsilon/2}$, the probability that $M_t$ hits $n^{0.5+\epsilon/2}$ before $-\frac{1}{2}n^{0.5+\epsilon}$ is at least 
\[ \frac{(\frac{p}{q})^{n^{0.5+\epsilon}-n^{0.5+\epsilon/2}}-1}{(\frac{p}{q})^{n^{0.5+\epsilon}-n^{0.5 + \epsilon/2}} - (\frac{p}{q})^{-2n^{0.5 + \epsilon/2}}} \geq 1 - e^{-\beta n^{1+2\epsilon}+2\beta n^{1 + 3\epsilon/2}}. \]
In particular, the probability of failure is exponentially small in $n$. By symmetry, the analogous statement holds for $M_t \leq \frac{1}{2}n^{0.5+\epsilon}$ when the spins on $A$ are set to $-1$. Started from $n^{0.5+\epsilon/2}$, the probability that $M_t$ hits $-n^{0.5+\epsilon/2}$ before $\frac{1}{2}n^{0.5+\epsilon}$ is at least $1 - e^{-(1-o(1))\beta n^{1+2\epsilon}}$ as well. Thus, the number of iterations that the magnetization of the system oscillates between $-\sqrt{n}$ and $\sqrt{n}$ under the corruption of $A$ dominates a geometric random variable with $e^{(1-o(1))\beta n^{1+2\epsilon}}$ mean. 

To conclude it suffices to show that $M_0$ is between $-n^{0.5+\epsilon/2}$ and $n^{0.5+\epsilon/2}$ with high probability. Since each spin is independent and uniform, $M_0 \overset{d}{=} 2\mathrm{Bin}(n - |A|, \frac{1}{2}) - (n-|A|)$. An application of Hoeffding's Inequality implies 
\[ \Prob{|M_0| > n^{0.5+\epsilon/2}} = \Prob{|2\mathrm{Bin}(n - |A|, \frac{1}{2}) - (n-|A|)| > n^{0.5+\epsilon/2}} \leq e^{-n^{\epsilon}/2}. \]
\end{proof}

\section{Maximizing the Adversary's Effect}
In this section, we present a series of results about the optimization problem from the perspective of the adversary. Namely, given a budget of $m$ nodes, what strategy should the adversary implement to induce the maximum effect at time $t$. We start with a result about the behavior at high temperature. The intuition in this regime is that the Glauber dynamics updates are well approximated by linear functions, and thus the degree of the vertices is the dominant contributing factor. In the remainder of this section, for $v$ a vertex in a graph $G$, let $N(v)$ denote the neighbors of $v$ and $d(v) = |N(v)|$ denote the degree. 
\begin{prop}
    Let $G$ be a graph of maximum degree $d$ and $M(\epsilon)$ be the $\epsilon$-mixing time for the Glauber dynamics on $G$ with inverse temperature $\beta$. Let $\tilde X_t^A$ be an instance of the $A$-corrupted Glauber dynamics for $|A| = m$. For $\beta \leq \frac{1}{2dm}$, $\epsilon = \frac{1}{4n}$ and $t \geq M(\epsilon)$, $\E{\left\vert\{v: \tilde X_t^A(v) = +1\}\right\vert}$ is maximized when $A$ are the $m$ vertices with largest degree. 
\end{prop}
\begin{proof}
    Fix a vertex $v$. Note that the effect of vertices at distance $k$ from $v$ is of order $\beta^k$. In particular, if we condition on the second neighborhood of $v$, then $\mu(\sigma_v = +1) = \frac{1}{2} \pm O((\beta d)^2)$. Thus, if we impose that the spins on $A$ are fixed to $+1$, then 
    \begin{align*}
        \mu(\sigma_v = +1 | \sigma_A = +1) &=
        \frac{1}{2} + \beta |N(v) \cap A| + \beta^2|N(N(v)) \cap A| + \cdots
    \end{align*}
    By linearity of expectation, 
    \[ \E{|\{v: \sigma(v) = +1\}|\middle\vert \sigma_A = +1} = \sum_v \mu(\sigma_v = +1 | \sigma_A = +1) = \frac{n}{2} + \beta\sum_{a \in A} d(a) + O(m(\beta d)^2) \]
    where the last equality follows from the maximum degree assumption and summing the resulting geometric series. For our choice of parameters, this lives on the interval 
    \[ \E{|\{v: \sigma(v) = +1\}|\middle\vert \sigma_A = +1} \in \left[\frac{n}{2} + \beta\sum_{a \in A} d(a) - \frac{1}{4}, \frac{n}{2} + \beta\sum_{a \in A} d(a) + \frac{1}{4} \right]. \]

    For $t \geq M(\epsilon)$, we know that the total variation distance between $\mu$ and the distribution of $X_t^A$ is at most $\epsilon$. Since $\left\vert\{v: \tilde X_t^A(v) = +1\}\right\vert$ is a random variable bounded by $n$, we know that 
    \[ \abs{\E{\left\vert\{v: \tilde X_t^A(v) = +1\}\right\vert} - \E{\left\vert\{v: \sigma_v = +1 | \sigma_A = +1\}\right\vert} } \leq \epsilon n = \frac{1}{4}. \]
    By the computations above, the triangle inequality implies that
    \[ \abs{\E{\left\vert\{v: \tilde X_t^A(v) = +1\}\right\vert} - \left(\frac{n}{2} + \beta \sum_{a \in A} d(a)\right) } \leq \frac{1}{2}.  \]
    It follows from this that in order to maximize $\E{\left\vert\{v: \tilde X_t^A(v) = +1\}\right\vert}$ it suffices to maximize the integer $\sum_{a \in A} d(a)$. 
\end{proof}
\begin{remark}
    For graphs $G$ and temperatures $\beta$ for which the Glauber dynamics mixes rapidly in polynomial time, we have that $M(\epsilon)$ for this choice of epsilon is also polynomial in $n$. In particular, the best choice for the adversary is to choose the highest degree vertices, and this choice manifests itself quickly. 
\end{remark}

In situations where the explicit optimizer is not as easy to compute, submodularity is a desirable property for the maximization of set functions. It is analogous to convexity, and thus allows for  efficient algorithms for approximate optimization, see for example the work of Feige, Mirrokni, and Vondr{\'a}k~\cite{feige2011maximizing}. It was shown by Bresler, Koehler, and Moitra~\cite{BrKoMo:19} using the GHS inequality that for the ferromagnetic Ising model, the influence maximization problem for the equilibrium measure is submodular. It is natural to ask whether or not this property extends to the dynamic, polynomial-time setting as well. The  following observations and results begin to answer this question. 
\begin{definition}[Submodularity]
    We say that a function $f: V \to \mathbb{R}$ is submodular if for all $S, T \subseteq V$, $f(S) + f(T) \geq f(S \cup T) + f(S \cap T)$.
\end{definition}

\begin{observation}
    The function $f_v: 2^{N(v)} \to [0,1]$ defined by $f_v(S) = \frac{e^{2|S| - d(v)}}{e^{2|S| - d(v)} + e^{-2|S| + d(v)}}$ is not submodular.
\end{observation}
The function $f_v$ is the local update of the Glauber dynamics at vertex $v$ when the set of neighbors of $v$ with $+1$ spin is exactly $S$. This function fails to be submodular when the sets $S$ and $T$ are too small. Thus, it immediately follows by considering a 10-regular graph and disjoint $S$ and $T$ of size 2 that the function $\E{\left\vert\{v: \tilde X_t^S(v) = +1\}\right\vert}$ cannot be a submodular function of $S$ in general for all graph $G$ and all times $t$. However, we know that as $t \to \infty$ the expectation approaches a submodular function. Thus, it is natural to ask whether or not there is a characterization for which graphs $G$, or for which times $t$, that the property of submodularity holds. We provide answers in two situations. 

\subsection{Approximate Submodularity}
\begin{definition}[Approximate Submodularity]
    We say that a function $f: V \to \mathbb{R}$ is \textit{$\epsilon$-approximately submodular} if for all $S, T \subseteq V$, $f(S) + f(T) \geq f(S \cup T) + f(S \cap T) - \epsilon$.
\end{definition}
\begin{prop}
    Let $M(\epsilon)$ be the $\epsilon$-mixing time of the Glauber dynamics on a graph $G$ with parameter $\beta$. In particular, for $t \geq M(\epsilon)$ we have $d_{\mathrm{TV}}(\tilde \mu_t^S, \mu) < \epsilon$ where $\tilde \mu_t^S$ is the distribution of $\tilde X_t^S$ and $\mu$ is the Ising model Gibbs measure. Then for $t \geq M(\epsilon)$, $\E{\left\vert\{v: \tilde X_t^S(v) = +1\}\right\vert}$ is $4\epsilon n$-approximately submodular as a function of $S$.
\end{prop}
\begin{proof}
    Note that the quantity $\left\vert\{v: \tilde X_t^S(v) = +1\}\right\vert$ is a random variable bounded below by 0 and above by $n$. The event where $X_t^S$ differs from $\mu$ has measure at most $\epsilon$, and the random variable differs by at most $n$ on this event, so the expectation under the two measures differs by at most $\epsilon n$. Applying this discrepancy to each term in the definition of submodularity yields the desired inequality. 
\end{proof}
Note once again that the above proposition in particular implies for rapidly mixing chains, the maximization problem becomes approximately submodular within a polynomial in $n$ number of time steps. 

\subsection{Submodularity Under Strong External Field}
We conclude by showing that in the presence of a strong enough external field, the maximization problem is in fact submodular for all time steps. The proof is an extension of local submodularity to global submodularity, using a coupling closely related to one introduced by Mossel and Roch~\cite{MosselRoch:10}. The following lemma establishes that the local update becomes submodular with a strong external field.
\begin{lemma}
    Let $h$ be at least the maximum degree of $G$. Let $f_v: 2^{N(v)} \to [0,1]$ be defined by $S \mapsto \frac{\exp(\beta(2|S|-d_v + h))}{\exp(\beta(2|S|-d_v + h)) + \exp(\beta(-2|S|+d_v - h))}$. Then $f_v$ is submodular.
\end{lemma}
\begin{proof}
    Note that $f_v(S) = \frac{1}{2}\left[ 1 + \tanh(\beta(2|S| - d_v + h)) \right]$. Thus, for submodularity, it suffices to prove that 
    \[ \tanh(\beta(2|S| - d_v + h)) + \tanh(\beta(2|T| - d_v + h)) \geq \tanh(\beta(2|S \cup T| - d_v + h)) + \tanh(\beta(2|S \cap T| - d_v + h)). \]
    We can compute that
    \begin{equation}\label{eq: tanh}
        \tanh(x) + \tanh(y) = \frac{e^x - e^{-x}}{e^x + e^{-x}} + \frac{e^y - e^{-y}}{e^y + e^{-y}} = 2\cdot \frac{e^{x+y} - e^{-x-y}}{e^{x+y} + e^{x-y} + e^{-x+y} + e^{-x-y}}.
    \end{equation}
    It holds that $|S| + |T| = |S \cup T| + |S \cap T|$ so we only need to compare $\exp(|S| - |T|) + \exp(-|S| + |T|)$ with $\exp(|S \cup T| - |S \cap T|) + \exp(-|S \cup T| + |S \cap T|)$. The function $e^x + e^{-x}$ is increasing for $x > 0$, and $|S| - |T| \leq |S \cup T| - |S \cap T|$ so 
    \[ \exp(|S| - |T|) + \exp(-|S| + |T|) \leq \exp(|S \cup T| - |S \cap T|) + \exp(-|S \cup T| + |S \cap T|). \]
    Substituting this into (1) and noting that $2|S|  + 2|T| - 2d_v + 2h \geq 0$ implies the desired inequality. 
\end{proof}

\begin{theorem}\label{thm: dynamic submodular}
    For any $S \subseteq V$, let $\tilde X_t^S$ be the $S$-corrupted Glauber dynamics with external field $h$ at least the maximum degree of $G$. The function $\E{\left\vert\{v: \tilde X_t^S(v) = +1\}\right\vert}$ is submodular function of $S$ for every $t$.
\end{theorem}
\begin{proof}
    Let $A_t = \{v \in V: \sigma_t^{S \cap T}(v) = +1\}$, $B_t = \{v \in V: \sigma_t^{S}(v) = +1\}$, $C_t = \{v \in V: \sigma_t^{T}(v) = +1\}$, and $D_t = \{v \in V: \sigma_t^{S \cup T}(v) = +1\}$. We construct a coupling such that $A_t \subseteq B_t \cap C_t$ and $D_t \subseteq B_t \cup C_t$.
    
    At each time step $t$, draw $v_t \sim \mathrm{Unif}(V \setminus A)$ and $\theta_t \sim \mathrm{Unif}([0,1])$. We define the following notation for the set of neighbors of $v_t$ with $+1$ spin in each chain:
    \begin{align*}
        A(v_t) &:= N(v_t) \cap (A_{t-1} \cup (S \cap T)) \\
        B(v_t) &:= N(v_t) \cap (B_{t-1} \cup S) \\
        C(v_t) &:= N(v_t) \cap (C_{t-1} \cup T) \\
        D(v_t) &:= N(v_t) \cap (D_{t-1} \cup S \cup T)
    \end{align*}
    Apply the following update rules to the respective chains:
    \[ \begin{cases}
        \sigma_t^{S \cap T}(v_t) = +1 & \text{if } \theta_t < f_{v_t}(A(v_t)) \\
        \sigma_t^{S}(v_t) = +1 & \text{if } \theta_t < f_{v_t}(B(v_t)) \\
        \sigma_t^{T}(v_t) = +1 & \text{if } \theta_t < f_{v_t}(A(v_t)) \text{ or } 1 - \theta_t < f_{v_t}(C(v_t)) - f_{v_t}(A(v_t)) \\
        \sigma_t^{S \cup T}(v_t) = +1 & \text{if } \theta_t < f_{v_t}(B(v_t)) \text{ or } 1 - \theta_t < f_{v_t}(D(v_t)) - f_{v_t}(B(v_t))
    \end{cases} \]
    and the spins are updated to $-1$ otherwise. Since $A_{t-1} \subseteq B_{t-1}$ by monotonicity $f_{v_t}(A(v_t)) \leq f_{v_t}(B(v_t))$, we have that $v_t \in A_t$ only if $v_t \in B_t$ and $v_t \in C_t$. Moreover if $v_t$ is removed from either $B_t$ or $C_t$ in this update, then it will be removed from $A_t$ as well. Thus, $A_t \subseteq B_t \cap C_t$. If $v_t \in D_t$, then either $\theta_t < f_{v_t}(B(v_t))$ in which case $v_t \in B_t$ or $1 - \theta_t < f_{v_t}(D(v_t)) - f_{v_t}(B(v_t))$, in which case $v_t \in C_t$ since 
    \[  f_{v_t}(D(v_t)) - f_{v_t}(B(v_t)) < f_{v_t}(C(v_t)) - f_{v_t}(A(v_t)). \]
    Indeed
    \begin{align*}
        f_{v_t}(B(v_t)) + f_{v_t}(C(v_t)) &\geq f_{v_t}(B(v_t) \cap C(v_t)) + f_{v_t}(B(v_t) \cup C(v_t)) \\
        &\geq f_{v_t}(A(v_t)) + f_{v_t}(D(v_t))
    \end{align*}
    where we used submodularity in the first inequality, and monotonicity and the induction hypothesis in the second inequality. Similarly if $v_t$ is removed from both $B_t$ and $C_t$ then it will be removed from $D_t$ as well, so $D_t \subseteq B_t \cup C_t$. 

    Now, 
    \begin{align*}
        |B_t| + |C_t| = |B_t \cap C_t| + |B_t \cup C_t| \geq |A_t| + |D_t|
    \end{align*}
    Taking expectations shows that $\E{\left\vert\{v: \tilde X_t^S(v) = +1\}\right\vert}$ is submodular. 
\end{proof}
\begin{remark}
    The above proof in general shows that for any submodular local update $f_v$, running the Glauber dynamics will result in a submodular optimization problem for all times $t$. 
\end{remark}

\printbibliography

\end{document}